\documentclass[12pt]{amsart}

\usepackage{amsmath, amscd, amssymb,xypic, euler}
\usepackage[frame,cmtip,arrow,matrix,line,graph,curve]{xy}
\usepackage{graphpap, color, paralist}
\usepackage[mathscr]{eucal}
\usepackage[pdftex,colorlinks,backref=page,citecolor=blue]{hyperref}

\hypersetup{%
pdftitle={Mori's program for the moduli space of pointed stable rational 
curves},%
pdfauthor={Han-Bom Moon},%
pdfkeywords={moduli space, birational geometry},%
}

\newtheorem{theorem}{Theorem}[section]

\theoremstyle{definition}

\newtheorem{remark}[theorem]{Remark}

\newtheorem*{acknowledgement}{Acknowledgement}

\newcommand{\PP}{\mathbb{P}}

\newcommand{\cO}{\mathcal{O} }
\newcommand{\cA}{\mathcal{A} }

\newcommand{\proj}{\mathrm{Proj}\;}

\def\Mzn{\overline{M}_{0,n} }
\def\Mza{\overline{M}_{0,\cA} }

\def\Mzek{\overline{M}_{0,n\cdot \epsilon_k} }

\def\git{/\!/ }

\begin{document}

\title
{Mori's program for the moduli space of pointed stable rational curves}
\date{November, 2010}
\author{Han-Bom Moon}
\address{Department of Mathematics, Seoul National University, Seoul 151-747, Korea}
\email{spring-1@snu.ac.kr}

\begin{abstract}
We prove that, assuming the F-conjecture, the 
log canonical model of the pair $(\Mzn, \sum a_{i}\psi_{i})$
is the Hassett's moduli space $\Mza$ without any modification of 
weight coefficients.
For the boundary weight cases, we prove that the birational model is 
the GIT quotient of the product of the projective lines.
This is a generalization of Simpson's theorem for symmetric weight cases.
\end{abstract}

\maketitle


\section{Introduction}\label{sec-introduction}
During the last several decades, the \emph{log canonical model} 
of a pair $(X, D)$ which is defined by 
\begin{equation}\label{eqn-logcanonicalmodel}
	X(K_{X}+D) := 
	\proj \left( \oplus_{l \ge 0} H^{0}(X, \cO(l(K_{X}+D)))\right)
\end{equation}
where the sum is taken over $l$ sufficiently divisible, plays
an important role in birational geometry especially in Mori's program.
In this paper, we prove the following theorem.
\begin{theorem}\label{thm-maintheorem}
Let $\cA = (a_{1},a_{2},\cdots, a_{n})$ be a weight datum, 
that is, a sequence of rational numbers such that $0 < a_{i}\le 1$.
Let $\Mzn$ be the moduli space of pointed stable rational curves.
Let $\psi_{i}$ be the $i$-th psi class (\cite[Section 2]{FarkasGibney}).
\begin{enumerate}
\item Assume the F-conjecture. Suppose that $\sum_{i=1}^{n} a_{i} > 2$.
Then $\Mzn(K_{\Mzn}+\sum_{i=1}^{n}a_{i}\psi_{i})$
is isomorphic to $\Mza$, the Hassett's moduli space of weighted pointed 
stable rational curves.
\item Suppose that $\sum_{i=1}^{n} a_{i} = 2$.
Then $\Mzn(K_{\Mzn}+\sum_{i=1}^{n}a_{i}\psi_{i})$ 
is isomorphic to the geometric 
invariant theory (GIT) quotient $(\PP^{1})^{n}\git_{L}SL(2)$ 
with respect to the linearization $L = \cO(a_{1}, \cdots, a_{n})$.
\end{enumerate}
\end{theorem}

For the definition and properties of $\Mza$, see \cite{Hassett}. 
The precise statement of the F-conjecture is in \cite[Question 1.1]
{KeelMcKernan}.
Note that item (2) does not rely on the F-conjecture.

Theorem \ref{thm-maintheorem} is an outcome of 
an attempt to generalize the following theorem of M. Simpson.
Set $m = \lfloor \frac{n}{2} \rfloor$. Let $\epsilon_{k}$ be a 
rational number in the range $\frac{1}{m+1-k} < \epsilon_{k} \le 
\frac{1}{m-k}$ for $k = 1, 2, \cdots, m-2$.
Let $n\cdot \epsilon = (\epsilon, \cdots, \epsilon)$ be a \emph{symmetric}
weight datum. 

\begin{theorem}\cite{AlexeevSwinarski, FedorchukSmyth, KiemMoon, Simpson}\label{thm-symmetriccase}
Let $\beta$ be a rational number satisfying $\frac{2}{n-1} < \beta
\le 1$ and let $D=\Mzn-M_{0,n}$ denote the total boundary divisor. 
\begin{enumerate} 
	\item If $\frac{2}{m-k+2} < \beta \le \frac{2}{m-k+1}$ for 
	$k = 1, 2, \cdots, m-2$, then $\Mzn(K_{\Mzn}+\beta D) \cong \Mzek$.
	\item  If $\frac{2}{n-1} < \beta \le \frac{2}{m+1}$, 
	then $\Mzn(K_{\Mzn}+\beta D) \cong (\PP^1)^n \git_{L} SL(2)$ 
	where $L = \cO(1,\cdots,1)$.
\end{enumerate}
\end{theorem}
Indeed, Theorem \ref{thm-symmetriccase} is a special case of 
Theorem \ref{thm-maintheorem} under the F-conjecture. See Remark
\ref{rmk-symmetriccase}.

In \cite{Fedorchuk}, Fedorchuk proved that for every weight 
datum $\cA$ and genus $g$, 
there exists a divisor $D_{g,\cA}$ on $\overline{M}_{g,n}$ such that 
$(\overline{M}_{g,n},D_{g,\cA})$ is a \emph{log canonical pair} and 
$\overline{M}_{g,n}(K_{\Mzn}+D_{g,\cA}) \cong \overline{M}_{g,\cA}$. 
His divisor $K_{\Mzn}+D_{0,\cA}$ is not proportional to 
$K_{\Mzn}+\sum_{i=1}^{n}a_{i}\psi_{i}$.

\begin{acknowledgement}
It is a great pleasure to thank my advisor Young-Hoon Kiem. 
It was impossible to finish this project without his patience and 
invaluable advice.
I would also like to thank Maksym Fedorchuk and David Swinarski 
who kindly answer to my questions.
\end{acknowledgement}


\section{Proof of the Main Theorem}\label{sec-proofofmainthm}

Throughout this section, we will assume $n \ge 4$. 
Fix a weight datum $\cA = (a_{1},a_{2},\cdots, a_{n})$
such that $\sum a_{i} > 2$.
Let $T'$ be the set of nonempty proper subsets of $[n] := \{1,2,\cdots,n\}$.
For each $I \in T'$ we can define the \emph{weight} of $I$ as
$w_{I} := \sum_{i \in I}a_{i}$.
Let $T \subset T'$ be the subset of $I \subset [n]$ such that 
$w_{I} < w_{I^{c}}$ or $w_{I} = w_{I^{c}}$, $|I| < |I^{c}|$ or 
$w_{I} = w_{I^{c}}$, $|I| = |I^{c}|$, $1 \in I$.
So there is a bijection between 
the set $\{ I \in T \;|\; 2 \le |I| \le n-2\}$ and 
the set of irreducible components $D_{I}$ of the boundary divisor of $\Mzn$.
Set $D_{I} = 0$ in the Neron-Severi vector space $N^{1}(\Mzn)$ 
for $I \in T$ with $|I| = 1$.

Let $\Delta_{\cA} := K_{\Mzn}+\sum_{i=1}^{n}a_{i}\psi_{i}$.
By \cite[Lemma 1]{FarkasGibney} and 
\cite[Lemma 3.5]{KeelMcKernan}, it is straightforward to show that 
\begin{equation}\label{eqn-divisor}
	\Delta_{\cA} = \sum_{\substack{I \in T \\ 2 \le |I| \le n-2}}
	\left(\frac{|I|(n-|I|)}{n-1}-2+
	\frac{(n-|I|)(n-|I|-1)}{(n-1)(n-2)}w_{I}+
	\frac{|I|(|I|-1)}{(n-1)(n-2)}w_{I^{c}}\right)D_{I}.
\end{equation}
By the geometry of reduction morphism $\varphi_{\cA} : \Mzn \to \Mza$
(\cite[\S 4]{Hassett}),
the following formulas are results of routine calculations.
Let $C \subset T$ (chosen for \emph{contracted}) 
be a subset such that $w_{I} \le 1$.
Let $D_{I} \in N^{1}(\Mza)$ be $\varphi_{\cA *}(D_{I})$
with abuse of notation.
\begin{equation}\label{eqn-pushforward}
	\varphi_{\cA *}(\Delta_{\cA}) = 
	\sum_{\substack{I \in C^{c}}}
	\left(\frac{|I|(n-|I|)}{n-1}-2+
	\frac{(n-|I|)(n-|I|-1)}{(n-1)(n-2)}w_{I}+
	\frac{|I|(|I|-1)}{(n-1)(n-2)}w_{I^{c}}\right)D_{I},
\end{equation}
\begin{multline}\label{eqn-pushforwardpullback}
	\varphi_{\cA}^{*}\varphi_{\cA *}(\Delta_{\cA})
	= \sum_{\substack{I \in C^{c}}}
	\left(\frac{|I|(n-|I|)}{n-1}-2+
	\frac{(n-|I|)(n-|I|-1)}{(n-1)(n-2)}w_{I}+
	\frac{|I|(|I|-1)}{(n-1)(n-2)}w_{I^{c}}\right)D_{I}\\
	+ \sum_{I \in C} \sum_{\substack{J \subset I\\|J| = 2}}
	\left(\frac{2(n-2)}{n-1}-2+\frac{(n-2)(n-3)}{(n-1)(n-2)}w_{J}
	+\frac{2\cdot 1}{(n-1)(n-2)}w_{J^{c}}\right)D_{I}\\
	= \sum_{\substack{I \in C^{c}}}
	\left(\frac{|I|(n-|I|)}{n-1}-2+
	\frac{(n-|I|)(n-|I|-1)}{(n-1)(n-2)}w_{I}+
	\frac{|I|(|I|-1)}{(n-1)(n-2)}w_{I^{c}}\right)D_{I}\\
	+ \sum_{I \in C}\left({|I| \choose 2}\left(\frac{2(n-2)}{n-1} - 2\right)
	+ \frac{n-3}{n-1}(|I|-1)w_{I}
	+ \frac{(|I|-1)(|I|-2)}{(n-1)(n-2)} w_{I}
	+ \frac{|I|(|I|-1)}{(n-1)(n-2)} w_{I^{c}}\right)D_{I}.
\end{multline}
Observe that if $|I|=1$ (so $I \in C$), 
then the coefficient of $D_{I}$ is zero.

From a direct calculation, it is immediate to check that 
\begin{equation}\label{eqn-difference}
	\Delta_{\cA} - \varphi_{\cA}^{*}\varphi_{\cA *}(\Delta_{\cA})
	= \sum_{\substack{I \in C\\2 \le |I| \le n-2}}
	\Big( (|I|-2)(1-w_{I})\Big) D_{I}.
\end{equation}
Note that for every $I \in C$ with nonzero $D_{I}$, 
$w_{I} \le 1$ by the definition of $C$.
So the difference $\Delta_{\cA} - 
\varphi_{\cA}^{*}\varphi_{\cA *}(\Delta_{\cA})$
is effective and supported on the exceptional locus of $\varphi_{\cA}$.
This implies that $	H^{0}(\Mzn, \Delta_{\cA}) \cong 
H^{0}(\Mzn, \varphi_{\cA}^{*}\varphi_{\cA *}(\Delta_{\cA})) \cong
H^{0}(\Mza, \varphi_{\cA *}(\Delta_{\cA}))$.
The same statement holds for a positive multiple of the same divisors.
Thus we get
\begin{equation}\label{eqn-lcmodeldescent}
	\Mzn(\Delta_{\cA})
	= \proj\left(\bigoplus_{l \ge 0} 
	H^0(\Mzn, \cO(l \Delta_{\cA}))\right)
	= \proj \left(\bigoplus_{l \ge 0}
	H^{0}(\Mza, \cO(l \varphi_{\cA *}(\Delta_{\cA})))\right).
\end{equation}
If we prove $\varphi_{\cA *}(\Delta_{\cA})$ is ample, 
then the last birational model is exactly $\Mza$.

Due to the F-conjecture (\cite[Question 1.1]{KeelMcKernan}), 
we may assume that the Mori cone $\overline{NE}_{1}(\Mzn)$
is generated by vital curve classes. 
This implies that $\overline{NE}_{1}(\Mza)$ is \emph{finitely} generated by 
the images of non-contracted vital curve classes.
Thus by Kleiman's criterion, a Cartier divisor $\varphi_{\cA *}(\Delta_{\cA})$ is 
ample on $\Mza$ if and only if for every non-contracted vital curve class 
$C \in \overline{NE}_{1}(\Mza)$, 
$C \cdot \varphi_{\cA *}(\Delta_{\cA}) > 0$.
By projection formula, this is equivalent to 
$\varphi_{\cA}^{*}\varphi_{\cA *}(\Delta_{\cA})$ is nef 
and it contracts exceptional curves only.
Therefore to check the ampleness of $\varphi_{\cA *}(\Delta_{\cA})$, 
it suffices to show that 
$\varphi_{\cA}^{*}\varphi_{\cA *}(\Delta_{\cA})$ intersects
positively with non-contracted vital curves.

Although the number of vital curve classes increases exponentially,
the intersection numbers of vital curves and 
$\varphi_{\cA}^{*}\varphi_{\cA *}(\Delta_{\cA})$ have surprisingly 
simple patterns.
For a partition $S_{1} \sqcup S_{2} \sqcup S_{3} \sqcup S_{4} = [n]$,
let $C(S_{1},S_{2},S_{3},S_{4})$ be the corresponding vital curve class.
Let $w_{j} := \sum_{i \in S_{j}}a_{i}$, $T_{j} := |S_{j}|$
and $w := w_{1}+w_{2}+w_{3}+w_{4}$.
We may assume that $w_{1} \le w_{2} \le w_{3} \le w_{4}$.

We encode the data of weights and partitions into 
a sequence of symbols of length $7$. 
Let $(s_{1},s_{2},s_{3},s_{4},s_{12},s_{13},s_{14})$ 
be a sequence of symbols such that each $s_{i}$ or $s_{ij}$ 
is one of $-, +, *$.
A vital curve class $C(S_{1},S_{2},S_{3},S_{4})$ is called 
of \emph{type} $(s_{1},s_{2},s_{3},s_{4},s_{12},s_{13},s_{14})$ if
the following conditions are satisfied:
\begin{compactitem}
	\item If $s_{i}$ (resp. $s_{ij}$) is $-$, then 
	$w_{i} \le 1$ (resp. $w_{i}+w_{j} \le 1$);
	\item If $s_{i}$ (resp. $s_{ij}$) is $+$, then $1 < w_{i} < w-1$
	(resp. $1<w_{i}+w_{j} < w-1$);
	\item If $s_{i}$ (resp. $s_{ij}$) is $*$, then 
	$w_{i} \ge w-1$ (resp. $w_{i}+w_{j} \ge w-1$).
\end{compactitem}

It is straightforward to see that a vital curve class 
$C(S_{1},S_{2},S_{3},S_{4})$ is 
contracted by $\varphi_{\cA}$ if and only if there exists $j$ such that 
$\sum_{i \notin S_{j}}a_{i} \le 1$.
So we may assume that $w_{i}+w_{j}+w_{k} > 1$ for all 
$\{i, j, k\} \subset \{1,2,3,4\}$, hence no $s_{i}$ is $*$.
Then we can divide all non-contracted vital curve classes into $13$ types.
By using \cite[Lemma 4.3]{KeelMcKernan}, 
we conclude that the formulas of intersection numbers 
$C(S_{1}, S_{2}, S_{3},S_{4})\cdot 
\varphi_{\cA}^{*}\varphi_{\cA *}(\Delta_{\cA})$
are depend only on the type of vital curves.
It is easy to calculate 
the intersection numbers by using a computer algebra system.
The list of all intersection numbers are in Table \ref{tbl-intersectionnumber}
in page \pageref{tbl-intersectionnumber}.
Note that all terms are nonnegative and for each intersection number
the last term is positive.
This completes the proof of item (1) of Theorem \ref{thm-maintheorem}.

\begin{table}[!ht]
\begin{tabular}{|c|p{0.75\textwidth}|}
\hline
Types & \multicolumn{1}{c|}{Intersection numbers}\\
\hline
$(-,-,-,-,-,-,-)$ & $(n-4)w_{1}+T_{1}(w-2)$\\ \hline
$(-,-,-,-,-,-,+)$ & $(T_{1}-1)(w_{1}+w_{2}+w_{3}-1)
+(T_{2}+T_{3}-2)w_{1}+(T_{4}-1)(1-w_{4})+(w-2)$\\ \hline
$(-,-,-,-,-,+,+)$  & $(T_{2}-1)w_{1}+(T_{1}-1)w_{2}
+(T_{3}-1)(1-w_{3})+(T_{4}-1)(1-w_{4})+(w-2)$\\ \hline
$(-,-,-,-,+,+,+)$ & $(T_{1}-1)(1-w_{1})+(T_{2}-1)(1-w_{2})
+(T_{3}-1)(1-w_{3})+(T_{4}-1)(1-w_{4})+(w-2)$\\ \hline
$(-,-,-,+,-,-,+)$ & $(T_{2}+T_{3}-2)w_{1}+
T_{1}(w_{1}+w_{2}+w_{3}-1)$ \\ \hline
$(-,-,-,+,-,+,+)$ & $(T_{2}-1)w_{1}+(T_{1}-1)w_{2}+
(T_{3}-1)(1-w_{3})+(w_{1}+w_{2}+w_{3}-1)$\\ \hline
$(-,-,-,+,+,+,+)$ & $(T_{1}-1)(1-w_{1})+(T_{2}-1)(1-w_{2})
+(T_{3}-1)(1-w_{3})+(w_{1}+w_{2}+w_{3}-1)$\\ \hline
$(-,-,+,+,-,+,+)$ & $T_{2}w_{1}+T_{1}w_{2}$\\ \hline
$(-,-,+,+,+,+,+)$ & $(T_{1}-1)(1-w_{1})+(T_{2}-1)(1-w_{2})+
w_{1}+w_{2}$\\ 
\hline
$(-,+,+,+,+,+,+)$ & $(T_{1}-1)(1-w_{1})+w_{1}+1$\\ \hline
$(+,+,+,+,+,+,+)$ & $2$\\ \hline
$(-,-,-,-,-,-,*)$ & $(T_{1}+T_{2}+T_{3}-3)(w_{1}+w_{2}+w_{3}-1)+
(T_{4}-1)(1-w_{4})+(w-2)$\\ \hline
$(-,-,-,+,-,-,*)$ & $(T_{1}+T_{2}+T_{3}-2)(w_{1}+w_{2}+w_{3}-1)$\\
\hline
\end{tabular}
\medskip
\caption{Intersection numbers for each vital curve types}
\label{tbl-intersectionnumber}
\end{table}

Next, we prove item (2) of Theorem \ref{thm-maintheorem}.
Set $T', T$ and $C$ as before. Define 
\begin{equation}\label{eqn-deltaprime}
	\Delta_{\cA}' := (n-4)\sum_{I \in T}
	\left(-{|I| \choose 2}\frac{2}{(n-1)(n-2)}+
	\frac{|I|-1}{n-2}w_{I}\right)D_{I}.
\end{equation}
Then it is straightforward to check that $\Delta_{\cA}-\Delta_{\cA}'$ is 
equal to the right side of \eqref{eqn-difference}.

By \cite{Kapranov}, there exists a birational morphism
$\pi_{\cA} : \Mzn \to (\PP^{1})^{n} \git_{L} SL(2)$ for 
any ample linearization $L = \cO(a_{1}, \cdots, a_{n})$.
Every boundary divisors except $D_{I}$ for 
$|I| = 2$ are contracted by $\pi_{\cA}$. 
The coefficients of $D_{I}$ in $\Delta_{\cA}-\Delta_{\cA}'$ is 
nonnegative since $w_{I} \le 1$ and is zero when $|I| = 2$.
Thus $\Delta_{\cA} - \Delta_{\cA}'$ is also effective 
and supported on the exceptional locus of $\pi_{\cA}$.
Therefore, by the same argument of the proof of item (1), 
$\Mzn(\Delta_{\cA})\cong \Mzn(\Delta_{\cA}')$.

For a vital curve class $C(S_{1},S_{2},S_{3},S_{4})$, 
by \cite[Lemma 4.3]{KeelMcKernan},
\begin{equation}\label{eqn-deltaprimeintersection}
	\Delta_{\cA}' \cdot C(S_{1}, S_{2}, S_{3}, S_{4}) = 
	\begin{cases}0, & w_{4} \ge 1\\
	(n-4)(1-w_{4}), & w_{4} \le 1 \mbox{ and } w_{1}+w_{4} \ge 1\\
	(n-4)w_{1}, & w_{4} \le 1 \mbox{ and } w_{1} + w_{4} \le 1.\\
	\end{cases}
\end{equation}
These intersection numbers are propotional to that of 
$\pi_{\cA}^{*}(\cO(a_{1}, \cdots, a_{n})\git SL(2))$
in \cite[Lemma 2.2]{AlexeevSwinarski}.
Since $\overline{NE}_{1}(\Mzn)$ is generated by vital curves, 
$\Delta_{\cA}'$ is proportional to the pull-back of the \emph{ample} 
divisor $\cO(a_{1}, \cdots, a_{n}) \git SL(2)$ on 
$(\PP^{1})^{n} \git_{L} SL(2)$.
Therefore $\Mzn(\Delta_{\cA})\cong (\PP^{1})^{n}\git_{L}SL(2)$.

\begin{remark}\label{rmk-symmetriccase}
The total psi-class $\psi := \sum_{i} \psi_{i}$ is
$\psi = \sum_{j=2}^{\lfloor n/2 \rfloor}\sum_{|I| = j} 
\frac{j(n-j)}{n-1} D_{I} = K_{\Mzn}+2D$ 
by \cite[Lemma 1]{FarkasGibney} and 
\cite[Lemma 3.5]{KeelMcKernan}.
So for $\alpha > 0$,
$K_{\Mzn}+\alpha \psi = (1+\alpha)(K_{\Mzn}+
\frac{2\alpha}{1+\alpha} D)$.
Therefore if $a_{1} = \cdots = a_{n} = \alpha$ 
for some $2/n < \alpha \le 1$, then 
$\Mzn(K_{\Mzn}+\sum a_{i}\psi_{i}) = \Mzn(K_{\Mzn}+\alpha \psi)$
is equal to $\Mzn(K_{\Mzn}+\frac{2\alpha}{1+\alpha}D)$.
If we substitute $\beta = \frac{2\alpha}{1+\alpha}$, then we 
get item (1) of Theorem \ref{thm-symmetriccase}.
Similarly, we can prove that $\Mzn(K_{\Mzn}+\frac{2}{n/2+1}D) \cong 
(\PP^{1})^{n}\git_{L} SL(2)$ with $L = \cO(2/n, 2/n, \cdots, 2/n)$ 
which is proportial to $\cO(1, 1, \cdots, 1)$.
So we can recover item (2) of Theorem \ref{thm-symmetriccase}
except the range of bigness.
\end{remark}

\begin{remark}\label{rmk-localglobal}
Theorem \ref{thm-maintheorem} shows a mysterious duality.
For $(C, x_{1}, x_{2}, \cdots, x_{n}) \in \Mzn$,
by the definition, the log canonical model $C(\omega_{C}+\sum a_{i}x_{i})$ 
is an $\cA$-stable curve and it is 
$\varphi_{\cA}(C, x_{1}, x_{2}, \cdots, x_{n})$. 
The \emph{same weight datum} determines 
$\Mzn(K_{\Mzn}+\sum a_{i}\psi_{i})$ 
of the moduli space $\Mzn$ itself.
\end{remark}

\begin{remark}\label{rmk-diversityofmodel}
In \cite{Keel}, Keel proved that 
$\mathrm{dim} \;N^{1}(\Mzn) = 2^{n-1}-{n \choose 2} - 1$.
By Theorem \ref{thm-maintheorem}, if the F-conjecture is true, 
then the family $\{\Mza\}$ of birational models of $\Mzn$ 
are detected by only an $n$-dimensional subcone 
of the effective cone of $\Mzn$.
So we can expect that still there is a huge wild world of unknown birational 
models of $\Mzn$.
\end{remark}


\begin{remark}\label{rmk-withoutFconjecture}
Although there is a strong belief on the F-conjecture, it seems that 
the proof of the F-conjecture is far from our hands.
So it is necessary finding a proof of Theorem \ref{thm-maintheorem}
without relying on the F-conjecture. 
As in the proof, proving the ampleness of $\varphi_{\cA *}(\Delta_{\cA})$ is 
a crucial step. We can express $\varphi_{\cA *}(\Delta_{\cA})$ 
in terms of tautological divisors on the universal curve of $\Mza$.
The author is working on proving the ampleness by using 
the expression and the technique of Fedorchuk in \cite{Fedorchuk}.
\end{remark}


\bibliographystyle{alpha}

\end{document}